\documentclass{ifacconf} 

\usepackage{mymacros} 
\begin{document}
\begin{frontmatter}

\title{On the Optimal Communication Weights in Distributed Optimization Algorithms} 
\thanks[footnoteinfo]{S. Colla is a FRIA grantee of the Fonds de la Recherche
Scientifique – FNRS. J. M. Hendrickx is supported by the F.R.S.-FNRS via its research project KORNET and its Incentive Grant for Scientific Research (MIS) “\emph{Learning from Pairwise Comparisons}”.}

\author[SebJul]{Sebastien Colla} 
\author[SebJul]{Julien M. Hendrickx} 

\address[SebJul]{ICTEAM Institute, UCLouvain, 1348 Louvain-la-Neuve, Belgium. \\(e-mail: \{sebastien.colla, julien.hendrickx\}@uclouvain.be).}

\begin{abstract} % Abstract of 50--100 words
We establish that in distributed optimization, the prevalent strategy of minimizing the second-largest eigenvalue modulus (SLEM) of the averaging matrix for selecting communication weights, while optimal for existing theoretical performance bounds, is generally not optimal regarding the exact worst-case performance of the algorithms. This exact performance can be computed using the Performance Estimation Problem (PEP) approach. We thus rely on PEP to formulate an optimization problem that determines the optimal communication weights for a distributed optimization algorithm deployed on a specified undirected graph. Our results show that the optimal weights can outperform the weights minimizing the second-largest eigenvalue modulus (SLEM) of the averaging matrix. This suggests that the SLEM is not the best characterization of weighted network performance for decentralized optimization. Additionally, we explore and compare alternative heuristics for weight selection in distributed optimization.  \vspace{-3mm}
\end{abstract}


\end{frontmatter}
%===============================================================================

\section{Introduction} \vspace{-1mm}
We seek to identify the best communication weights to use in decentralized optimization algorithms, by leveraging the Performance Estimation Problem framework (PEP) from \cite{Colla_PEP_dec}.  \vspace{-1mm}

In decentralized optimization, one considers a set of agents $\V = \{1,\dots,n\}$, connected by a set $E$ of $m$ communication links via the network $G(\V,E)$.
The agents seek to collaboratively minimize the average of their private local functions $f_i: \mathbb{R}^{d}\to\mathbb{R}$: \vspace{-3mm}
\begin{align} \label{opt:dec_prob}
      \underset{\text{\normalsize $x \in\mathbb{R}^{d}$}}{\mathrm{minimize}} \quad f(x) = \frac{1}{n}\sum_{i=1}^n f_i(x). \vspace{-3.5mm}
\end{align}
Each agent $i$ holds a local copy $x_i$ of the decision variable to perform local computations. The agents exchange local information with their neighbors in $G$ to gradually come to an agreement on the minimizer $x^*$ of the global function $f$. These exchanges often take the form of an average consensus on some quantity, e.g.,  on the $x_i$. The consensus step can be represented using multiplication by an averaging matrix $W \in \mathbb{R}^{n\times n}$, for which $W_{ij} = 0$ if there is no communication link between $i$ and $j$. In this work, we focus on the case where the communication network is undirected, which corresponds to a symmetric matrix $W$. All the assumptions for $W$ are summarized below 
\begin{assumption}[Averaging matrix] 
   \label{as:mat}
   The averaging matrix
   $W = [W_{ij}] \in \Rmat{n}{n}$ satisfies \vspace{-1mm}
   \begin{enumerate}
      \item $W^T = W$,  \hfill  (Symmetry) 
      \item  $W \mathbf{1} = \mathbf{1}$ and $\mathbf{1}^T W  = \mathbf{1}^T$, \hfill (Averaging Consensus)
      \item $W \in \T$,
      \hfill (Topology) \\
      where $\T = \{W : W_{ij} = 0 \text{ if $(i,j) \notin E$ and $i \ne j$}\}$.
   \end{enumerate}
\end{assumption}
If, in addition to (2), the matrix $W$ is non-negative, we say that $W$ is doubly-stochastic. 
While it is a common assumption in the literature, we choose not to use it because most decentralized convergence results do in fact not use the non-negativity assumption. Moreover, non-negativity is not necessary for the convergence of a pure consensus protocol, as shown in \cite{xiaoboyd2004fast}.

An example of a well-known distributed optimization algorithm is DIGing, from \cite{DIGing}. The algorithm uses a gradient tracking variable $s_i^k$ and can be written, for each agent $i \in \V$, as
   \begin{align}
      x_{i}^{k+1} &= \sum_{j} W_{ij}\, x_{j}^k - \alpha s_i^k, \label{eq:DIGing1} \\
      s_{i}^{k+1} &= \sum_{j} W_{ij}\, s_{j}^{k} + \nabla f_i(x_i^{k+1}) - \nabla f_i(x_i^{k}), \label{eq:DIGing2}
   \end{align}
where $\alpha > 0$ is a constant step-size.

In general, the performance of a decentralized optimization method is largely impacted by the averaging matrix $W$, see for example the survey \cite{DGD}. While the zero elements are imposed by the network topology, the values of the non-zero elements can be freely chosen, and could be carefully determined to obtain efficient algorithms, as for any other parameter of a method.

Almost all theoretical performance guarantees from the literature of decentralized optimization algorithms depend on the second-largest singular value of $W$ (after one), denoted $\lam$, or some equivalent measure, because it is easy to use in the proofs and characterizes well the network behavior in terms of consensus.
By Assumption \ref{as:mat}, $W\mathbf{1}=\mathbf{1}$ and therefore, $\lam$ can be computed as
\begin{equation} \label{eq:lam} 
   \lam = \|W - \mathbf{11}^T/n \|_2.
\end{equation}
Decentralized algorithms require $\lam <1$, otherwise, the agents never converge to a common point, see \cite{xiaoboyd2004fast}. This requirement is not guaranteed by Assumption \ref{as:mat} but is implicit since we are looking for averaging matrices that enable fast convergence for distributed optimization algorithms, thus ruling out matrices that prevent convergence. 
When $\lam < 1$ and Assumption \ref{as:mat} hold, the eigenvalues of the symmetric matrix $W$ are as follows
$$ 1 = \lam_1(W) \le \lam_2(W) \le \dots \le \lam_n(W). $$
In this case, $\lam$ also corresponds to the Second-Largest Eigenvalue Modulus (SLEM) of $W$, 
\begin{equation} \label{eq:SLEM} 
   \lam = \max\bigl\{\,|\lam_2(W)|,~ |\lam_n(W)|\,\bigr\}.
\end{equation}

All known theoretical performance guarantees improve when $\lam$ decreases. Therefore, the classical theoretical approach for selecting the weights in distributed optimization is to choose those leading to the smallest possible $\lam$ for the given network topology, which results in the smallest error guarantee.

\subsection*{Contributions}
In this paper, we formulate an optimization problem that computes the optimal communication weights for a given distributed optimization algorithm, a given network topology, and other given settings. The formulation relies on the Performance Estimation Problem framework (PEP) \cite{PEP_Smooth} which allows to numerically compute the exact worst-case performance of an optimization algorithm, and which has been extended to distributed optimization in \cite{Colla_PEP_dec}. The resulting weight-tuning problem has no guarantee to be smooth nor convex, but we can obtain good estimates of the solution via a zero-order method. We show that in many settings, the optimal weights are different from those minimizing $\lam$ and allow, for example, to decrease the convergence time by up to 4 for the DIGing algorithm \cite{DIGing}. It seems that all the eigenvalues, as well as their sign, are important to determine the performance of a decentralized optimization method. This indicates that the second-largest eigenvalue modulus $\lam$ (SLEM) is not the best determinant of weighted network performance in distributed optimization methods. We also explore different weights heuristics and compare them with the optimal ones, to find a better characterization of network performance in distributed optimization methods. 

\section{Weights Heuristics} %Heuristics to choose the weights} 
\label{sec:heuristics}
We review possible heuristics to set the communication weights in distributed optimization, for which we will analyze the resulting performance in Section \ref{sec:results}.

Let $B \in \Rmat{n}{m}$ be the oriented incidence matrix of the network $G$, defined as
$$ B_{ie} = 
\begin{cases} 
   1  & \text{ if edge $e$ starts from node $i$,}\\
   -1 & \text{ if edge $e$ ends at node $i$,} \\
   0 & \text{ otherwise.}\\
\end{cases}
$$
Since the edges are undirected, any orientation can be chosen for each edge. The matrix $B$ allows defining an averaging matrix $W$ satisfying Assumption \ref{as:mat}, for any choice of the edge weights
\begin{equation} \label{eq:Wdef}
    W = I - B \mathrm{diag}(w) B^T,
\end{equation}
where $w \in \Rvec{m}$ is the vector of weights for the $m$ edges in $E$ and $\mathrm{diag}(w) \in \Rmat{m}{m}$ is a diagonal matrix with $w$ on the diagonal. 

Common choices for the averaging matrix $W$ in decentralized optimization are those developed for the pure linear averaging consensus protocols, 
\begin{equation} \label{eq:cons}
   \mathbf{x}^{k+1} = W \mathbf{x}^k,
\end{equation}
and are listed below. \vspace{-3mm}

\paragraph*{Minimum-$\lam$ weights.} In \cite{xiaoboyd2004fast}, the authors have proved that the per-step convergence factor for the averaging consensus \eqref{eq:cons} is given by $\lam=\|W - \mathbf{11}^T/n \|_2$, and therefore, that the fastest averaging matrix for the consensus, denoted $\Wls$, is obtained by minimizing $\lam$, %defined in \eqref{eq:lam}:
\begin{align} \label{eq:minlam2}
   \Wls &= I - B \mathrm{diag}(w_{\lam^*}) B^T, \\
   \text{where} ~w_{\lam^*} &= \underset{w \in \Rvec{m}}{\mathrm{argmin}} \| I - B \mathrm{diag}(w) B^T - \mathbf{11}^T/n\|_2. \label{eq:minlam_w}
\end{align}
While computing $\Wls$ requires global knowledge of the network, some papers have explored ways to compute it in a decentralized way, e.g. \cite{boyd2006randomized}. \vspace{-3mm} 

\paragraph*{Uniform edge weights.} When we choose all the weights to be equal to $q \in \R$, then we have $w = q\mathbf{1}$ and
$$ \Wq{q} = I - q BB^T.$$
According to \cite{xiaoboyd2004fast}, the value of $q$ minimizing the convergence factor, given by $\lam$, is given by 
\begin{equation} \label{eq:constant-weight}
   q_{\lam^*} = \frac{2}{\lam_2(L) + \lam_{n}(L)},
\end{equation}
where $L = BB^T$ is the Laplacian matrix of the graph $G$ and its eigenvalues $\lam_i(L)$ are labelled in ascending order. When the network $G$ is connected, the resulting averaging matrix guarantees the convergence of the averaging consensus \eqref{eq:cons}, i.e. $\lam < 1$.\vspace{-3mm}
\paragraph*{Maximum-degree weights.} Another uniform edge weight that
always yields asymptotic average consensus is
\begin{equation}
   \label{eq:max-degree-weights}
   q_{d} = 1/(d_{\max}+1),
\end{equation}
where $d_{\max}$ is the maximal degree in $G$, which can easily be computed in a decentralized manner. 
\vspace{-3mm} 

\paragraph*{Metropolis weights.} The Metropolis weights only use local-degree information so that each node can set up its weights without knowing global information on the network. They are defined, for each edge $e = (i,j) \in E$, as
\begin{equation} \label{eq:metropolis}
    w_{e} = \frac{1}{\max\{d_i, d_j\}+1},
\end{equation}
where $d_i$ is the degree of node $i$, without counting the node itself. The averaging matrix calculated from \eqref{eq:Wdef} with Metropolis weights is denoted $W_M$ and guarantees the convergence of the averaging consensus \eqref{eq:cons}. These weights are derived from the Metropolis–Hastings algorithm, when employed for the simulation of a Markov chain with uniform equilibrium distribution, see \cite{boyd2004fastest}. \vspace{-3mm}
%% NOTE: IF the graph is NOT bipartite, then we can remove the '+1' after the degree.

\paragraph*{Lazy Metropolis weights.} A popular variation of Metropolis weights, introduced by \cite{olshevsky2015linear}, is defined as %is the \emph{lazy Metropolis weights},
$$ W_{\mathrm{lazy}-M} = \frac{1}{2} (I + W_M).$$
This averaging matrix is diagonally dominant and has thus non-negative eigenvalues, which can be useful for proofs.

 In the pure consensus \eqref{eq:cons}, the decrease of the consensus error $\|\x^k - \xb\mathbf{1}\|$ is dominated by $\lam$ because the influence of all the other smaller eigenvalues (in absolute value) vanishes more rapidly, and the influence of $\lam_1=1$ is zero because the error vector is orthogonal to the associated eigenvector $\mathbf{1}$. In decentralized optimization, the consensus is perturbed by local gradient updates, which prevent the vanishing of smaller eigenvalues of the averaging matrix. Therefore, all the eigenvalues may influence the convergence speed of decentralized optimization algorithms. This motivates us to consider weight heuristics based on all the eigenvalues of the averaging matrices, and that we consider potentially relevant in decentralized optimization. \vspace{-4.5mm}

\paragraph*{Minimum-$\|W\|_{\Sigma}$ weights.}
The nuclear norm of a symmetric matrix $W$ is defined as 
$$ \|W\|_{\Sigma} = \sum |\lam_i(W)|, $$
and the weighted nuclear norm 
is given by
$$ \|W\|_{\gamma,\Sigma} = \sum \gamma_i|\lam_i(W)|, $$
where $\gamma = [\gamma_1,\dots,\gamma_n]$, with $\gamma_i \ge 0$ a non-negative weight assigned to $\lambda_i$. Both norms are used as convex surrogates for the matrix rank in problems involving rank minimization. The weights minimizing $\|W\|_{\gamma,\Sigma}$ can be relevant in distributed optimization because it generalizes the weights minimizing $\lam$, by taking into account all the other eigenvalues, which may play a role in the algorithm convergence. We denote the corresponding averaging matrix as $\Wals$ and define it formally as
\begin{align}
   \Wals &= I - B \mathrm{diag}(\wals) B^T, \\
   \text{where} ~\wals &= \underset{w \in \Rvec{m}}{\mathrm{argmin}} \| I - B \mathrm{diag}(w) B^T\|_{\gamma,\Sigma},\\[-5mm]
\end{align}
which can be computed by solving an equivalent SDP, described in \cite{alizadeh1995interior}. \vspace{-4.5mm}
\paragraph*{Minimum-$\Rtot$ weights.} The total effective resistance of a weighted graph, denoted $\Rtot$, can be defined based on the eigenvalues of the associated averaging matrix $W$: \vspace{-1mm}
\begin{equation} \label{eq:Rtot} 
   \Rtot(W) = n \sum_{i=2}^n \frac{1}{1-\lam_i(W)}. 
\end{equation}
The initial interpretation and definition of $\Rtot$ is related to the total electrical resistance of a resistor network with conductances given by the edges. \cite{ghosh2008resistance} provide many useful interpretations of the total effective resistance. For example, $\Rtot$ is related to the average commute time in the Markov Chain derived from the weighted graph. 
Hence, the averaging matrix $\Wr$ minimizing the total effective resistance would consist in the graph leading to the smallest average commute time:
\begin{align} \label{eq:minRes}
   \Wr &= I - B \mathrm{diag}(w_{\Rtot^*}) B^T, \\
   \text{where}~w_{\Rtot^*} &= \underset{w \in \Rvec{m}}{\mathrm{argmin}} ~\Rtot(I - B \mathrm{diag}(w) B^T) \\
      &\text{s.t.} \sum_{e \text{ adjacent to } i} \hspace{-5mm} w_e \le 1 - \epsilon \quad \text{for all $i$},
\end{align}
where $\epsilon \in [0,1)$ is the smallest accepted value for the diagonal elements of the averaging matrix. Using results from \cite{ghosh2008resistance}, the solution to problem \eqref{eq:minRes} can be computed by solving an equivalent SDP. \vspace{-4.5mm}

\paragraph*{Least mean-square deviation weights.}
\cite{xiaoboyd2007distributed} formulate a convex problem to compute the symmetric averaging weights that minimize the mean-square deviation from the mean produced by the average consensus with additive noise. This extension of the averaging iteration \eqref{eq:cons} adds a noise at each node and each step: \vspace{-1mm}
\begin{equation} \label{eq:cons_noise}
   \mathbf{x}^{k+1} = W \mathbf{x}^k + \mathbf{v}^k,
\end{equation}
where $\mathbf{v}^k = [v_1^k \dots v_n^k]^T$, and $v_i^k$ are independent random variables, identically distributed, with zero mean and unit variance. Due to the noise, the agent values $x_i$ are not converging to the average, but the quality of the consensus can be evaluated using the mean-square deviation from the average $\xb = \mathbf{x}^T \mathbf{1}/n$: \vspace{-2mm}
$$\delta^k = \Exp \sum_{i=1}^n (x_i^k - \xb)^2 $$
According to \cite{xiaoboyd2007distributed}, $\delta^k$ converges to a steady state mean-square deviation which depends on the eigenvalues of $W$ as \vspace{-1mm}
\begin{equation} \label{eq:deltass}
    \delta_{ss}(W) = \lim_{k\to\infty} \delta^k = \sum_{i=2}^n \frac{1}{1-\lam_i(W)^2}.
\end{equation}
The averaging matrix $\Wlms$ minimizing the steady state mean-square deviation is therefore given by
\begin{align} \label{eq:least_dev_W}
   \Wlms &= I - B \mathrm{diag}(w_{\delta^*}) B^T, \\
   \text{where} ~w_{\delta^*} &= \underset{w \in \Rvec{m}}{\mathrm{argmin}} ~\delta_{ss}(I - B \mathrm{diag}(w) B^T) \label{eq:least_dev} 
\end{align}  
\cite{xiaoboyd2007distributed} show that problem \eqref{eq:least_dev} is convex and details how to evaluate the function $\delta_{ss}$, as well as its gradient and hessian, such that it can be minimized using standard optimization methods. Choosing $\Wlms$ for distributed optimization may be relevant since, distributed optimization algorithms are similar to the noisy consensus \eqref{eq:cons_noise}, where additive noise is replaced by a local update based on local gradients information, that may not be i.i.d., of zero mean or unit variance.

\section{Worst-case Performance evaluation} \label{sec:PEP} \vspace{-1.5mm}
To evaluate the quality of an averaging matrix $W$ for a given decentralized algorithm, we compute the worst-case performance it provides. 
For this purpose, we rely on the Performance Estimation Problem framework (PEP) that was initially developed for centralized optimization by \cite{PEP_Smooth} and extended to decentralized optimization by \cite{Colla_PEP_dec}.
To obtain a tight bound on the performance of an algorithm $A$, the conceptual idea is to search for instances of local functions and starting points for all agents, allowed by the setting considered, producing the largest error after a given number $K$ of iterations of the algorithm. This idea can be formulated as a real optimization problem that maximizes the error measure $\P$ of the algorithm result, over all possible functions and initial point allowed: 
\begin{align}
   \hspace{-1mm} E_{\S}(W,\alpha)=&\underset{x^*,\{x_i^0, f_i\}_{i\in \V}}{\max} ~ \P(f_i,x_i^0,\dots,x_i^K,x^*) & \label{eq:PEP}\\
   \text{ s.t. } \quad& x_i^{k} \quad \text{from algorithm A} & \hspace{-15mm}\text{$\substack{k=1,\dots,K \\ i \in \V\qquad}$ \small{~(algorithm)}} \\
   &  f_i \in \F,& \hspace{-25mm} \text{\small{ $i \in \V$ ~ (class of functions)}} \\%[1mm]
   & x_i^0 \quad \text{satisfies $\I$}& \hspace{-15mm}\text{\small{(initial condition)}} \\
   & \hspace{5.5mm} \text{\small{(e.g. $\|x_i^0 - x^*\| \le 1$,~ $i \in \V$)}}  & \\
   & \frac{1}{n} \sum_{i=1}^n \nabla f_i(x^*) = 0, & \hspace{-15mm} \text{\small{(optimality condition)}}
\end{align}
where $\alpha$ is the step-size parameter of the algorithm and $\S$ is the performance evaluation setting which specifies the graph of agents, the algorithm, the number of steps, the performance criterion, the class of function, and the initial conditions: \vspace{-2mm}
$$\S = (G,A,K,P,\F,\I).$$
The worst-case performance function $E_{\S}(W,\alpha)$ is always specific to a given setting $\S$, but it may sometimes be ignored in the notation $E(W,\alpha)$.
Problem \eqref{eq:PEP} can be solved exactly via an SDP reformulation, relying on a discretization of the functions and the use of interpolation constraints appropriate to the given function class $\F$, see \cite{PEP_composite}. We can use any of the common classes of functions such as $\mu$-strongly convex and $L$-smooth functions. Moreover, the performance measure $P$ and the initial conditions $\I$ can be any expression that is linear in the local function values and quadratic in the iterates and gradient values, see \cite{Colla_PEP_dec} for more details. \vspace{-0.5mm}

The PEP framework \eqref{eq:PEP} allows trying the different weight heuristics from Section \ref{sec:heuristics} on different graphs and algorithms, and to compare their resulting worst-case performance. Furthermore, it also allows for seeking the weights leading to the best performance.
 
\section{Weights minimizing the worst-case performance} \label{sec:tuning} \vspace{-0.5mm}
To find the best averaging matrix to use for a setting $\S$, we can find the matrix $W^*$, along with a step-size $\alpha^*$, minimizing the performance function $E_{\S}(W,\alpha)$ of the algorithm, defined in \eqref{eq:PEP},
\begin{align} \label{eq:tuning}
   (W^*, \alpha^*) = \underset{\substack{W \in \Rmat{n}{n}\\ \alpha \ge 0 \quad}}{\mathrm{argmin}} & ~ E_{\S}(W,\alpha) \\
   \text{s.t. }   
   & W = W^T, \quad  W\mathbf{1} = \mathbf{1}, \quad W \in \T. \\
   & \lam(W) = \| W - \mathbf{11}^T/n \|_2 < 1. \label{eq:lamconv}
\end{align}
Constraint $\lam(W) < 1$ is necessary for the convergence of decentralized algorithms. Matrices that do not satisfy it should, in principle, lead to larger values of $E_{\S}$, when the number $K$ of iterations is large enough. To ensure the condition is met, even for small $K$, we interpret $E_{\S}$ as $\infty$ when it does not hold. Therefore, Problem \eqref{eq:tuning} can be written without constraints, using \eqref{eq:Wdef},
\begin{align} \label{eq:tuning2}
   (w^*, \alpha^*) = \underset{\substack{w \in \Rvec{m}\\ \alpha \ge 0 \quad}}{\mathrm{argmin}} & ~ E_{\S}( I - B\mathrm{diag}(w)B^T,\alpha), 
\end{align}
with step-size $\alpha \in \R$ and edge weights vector $w \in \Rvec{m}$ as variables, and where we interpret $E_{\S}(W,\alpha)$ as $\infty$ whenever $\| W - \mathbf{11}^T/n \|_2 \ge 1$.
This problem is a priori non-convex and non-smooth. While this is difficult to prove, it could be smooth or convex for some settings.
A simple way of solving Problem \eqref{eq:tuning} is using a zero-order method, that only relies on function evaluations to identify a minimizer. In this work, we use the pattern search method from Matlab, which is described and analyzed in \cite{audet2002pattern}. Their analysis showed that even if the objective function is discontinuous or extended-valued, the method finds a limit point with some minimizing properties. Moreover, if the function is strictly differentiable at this limit point, then it corresponds to a local minimum. We have, of course, no guarantee to converge to a global minimizer in the non-convex case. There are other possible approaches to solve \eqref{eq:tuning}: \vspace{-0.5mm}
\begin{itemize}
   \item Problem \eqref{eq:tuning} corresponds to a min-max problem. The SDP reformulation of the inner maximization \eqref{eq:PEP} can be dualized, so that it can be combined with the outer minimization. The resulting problem is a quadratically constrained quadratic program (QCQP), which is known to be NP-hard in general. However, there are solvers, such
   as Gurobi, that find globally optimal solutions to non-convex QCQPs in finite time.
   \item One could try to solve Problem \eqref{eq:tuning} using first-order methods, which would allow converging to a local minimum. The generalized gradient of function $E(W,\alpha)$ can be computed via an SDP sensitivity analysis, see \cite{bonnans2013perturbation}. Indeed, as explained in \cite{Colla_PEP_dec}, Problem \eqref{eq:PEP} can be reformulated in an SDP.
\end{itemize}
For this first analysis, we have decided to use a simple zero-order method, which already provides interesting results, as shown in the following section. \vspace{-1mm}

\section{Results} \label{sec:results} \vspace{-1mm}
\begin{figure*}[t]
   \vspace{-1mm}
   \centering
   \begin{tabular}{C{2.5cm}|*{4}{C{3.4cm}}}
   
      \hspace*{-2.5mm}\includegraphics[width=2.9cm]{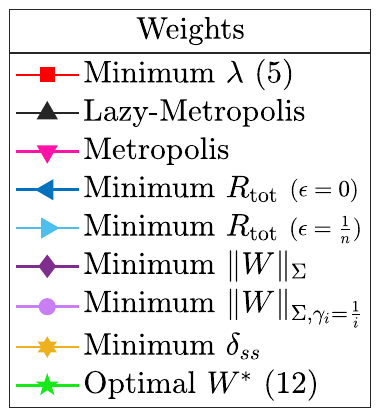} 
      & \vspace{-12mm} \begin{tabular}{@{}c@{}} Complete \\[2mm] \includegraphics[width=2cm]{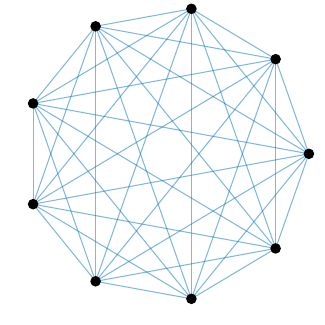} \\ ($n=9$) \end{tabular} 
      & \vspace{-12mm} \begin{tabular}{@{}c@{}} Star \\[2mm] \includegraphics[width=2cm]{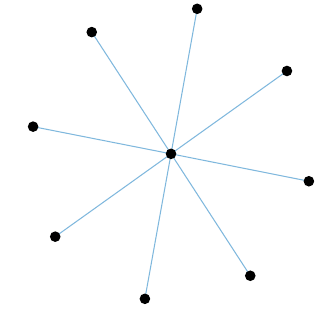} \\ ($n=9$) \end{tabular} 
      & \vspace{-12mm} \begin{tabular}{@{}c@{}} Cycle \\[2mm] \includegraphics[width=2cm]{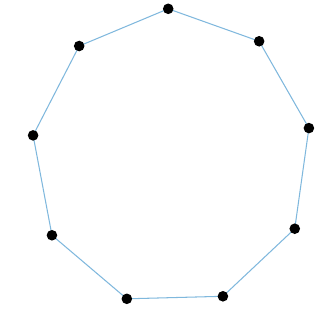} \\ ($n=9$) \end{tabular} 
      & \vspace{-12mm} \begin{tabular}{@{}c@{}} Grid \\[2mm] \includegraphics[angle=3,origin=c,width=2cm]{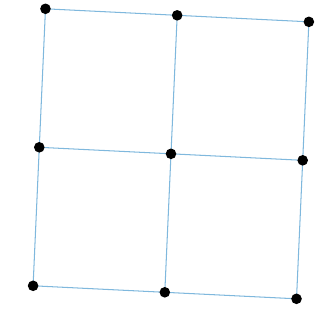} \\ ($n=9$) \end{tabular} \\ \hline

      \begin{tabular}{@{}l@{}}
         DIGing \\ 
         \scriptsize{$E_{\S} = $ conv. rate on} \\
         \scriptsize{$\frac{1}{n} \sum_{i=1}^n \|x_i^k-x^*\|^2$}
      \end{tabular} 
     & \includegraphics[width=3.5cm]{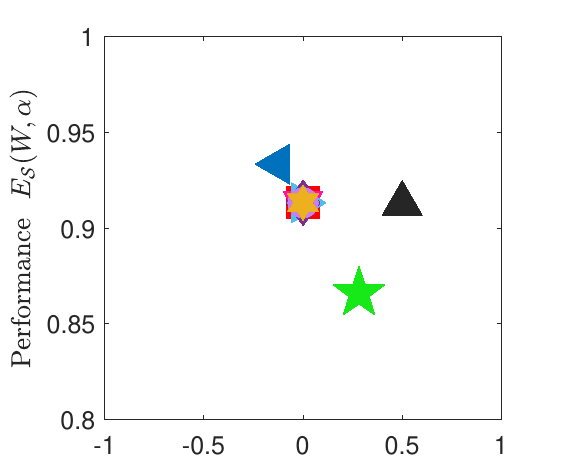} 
     & \includegraphics[width=3.5cm]{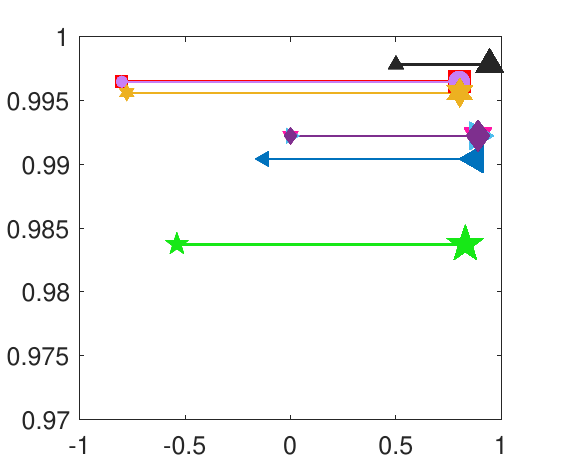} 
     & \includegraphics[width=3.5cm]{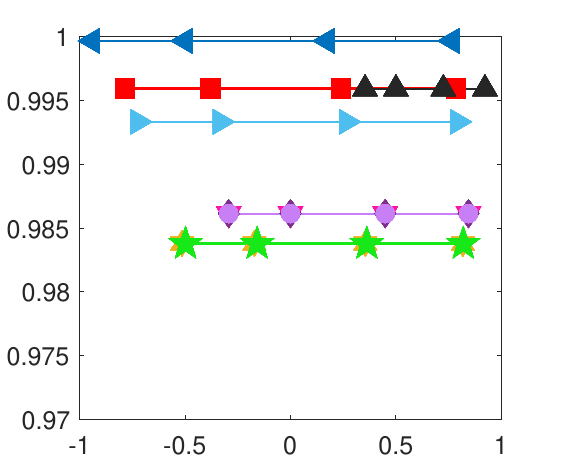} 
     & \includegraphics[width=3.5cm]{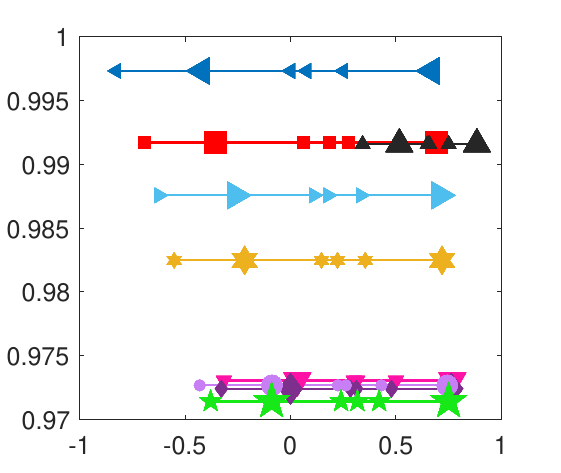} \\

     \begin{tabular}{@{}l@{}}
      ATC-DIGing  \\
      \scriptsize{$E_{\S} = $ conv. rate on} \\
      \scriptsize{$\frac{1}{n} \sum_{i=1}^n \|x_i^k-x^*\|^2$}
   \end{tabular}
      & \includegraphics[width=3.5cm]{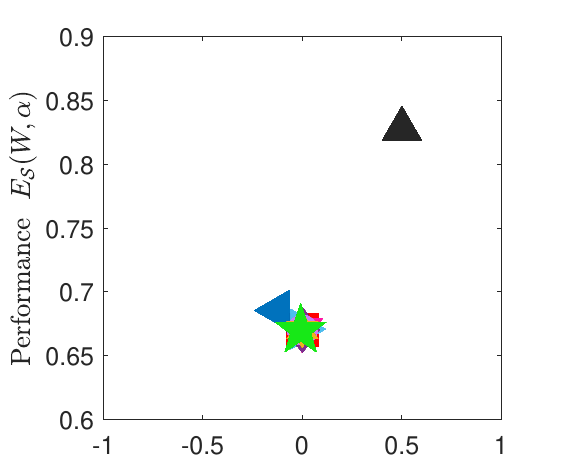} 
      & \includegraphics[width=3.5cm]{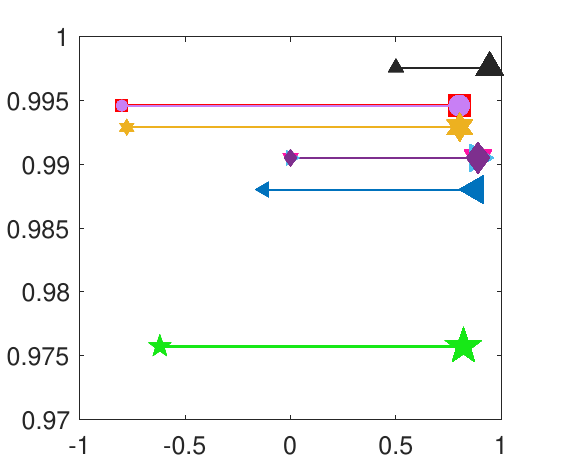} 
      & \includegraphics[width=3.5cm]{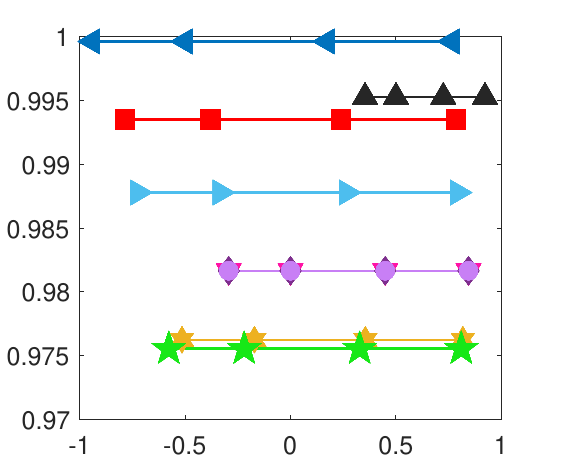} 
      & \includegraphics[width=3.5cm]{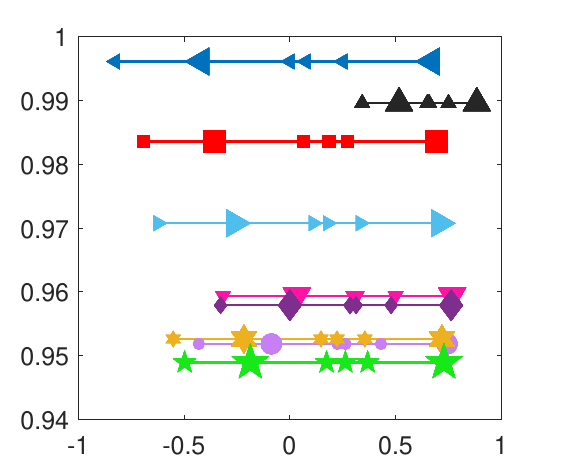} \\
 
      \begin{tabular}{@{}l@{}}
         EXTRA \\ 
         \scriptsize{$E_{\S} = $ worst-case of} \\
         \scriptsize{$f(\xb^K)-f(x^*)$} \\
         \scriptsize{for $K=5$.}
      \end{tabular}
      & \includegraphics[width=3.5cm]{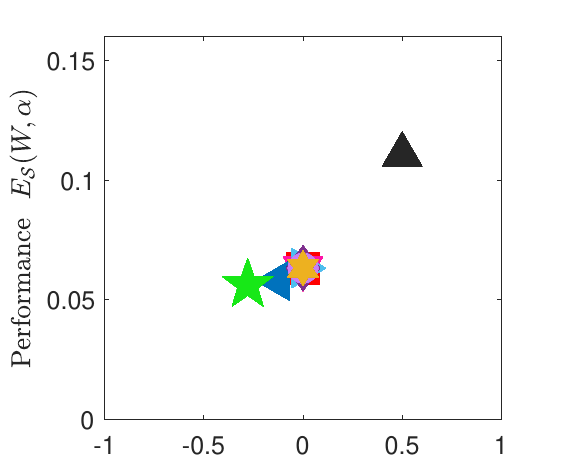} 
      & \includegraphics[width=3.5cm]{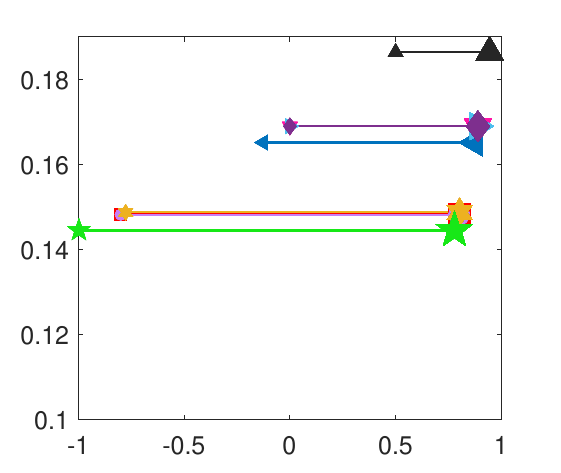} 
      & \includegraphics[width=3.5cm]{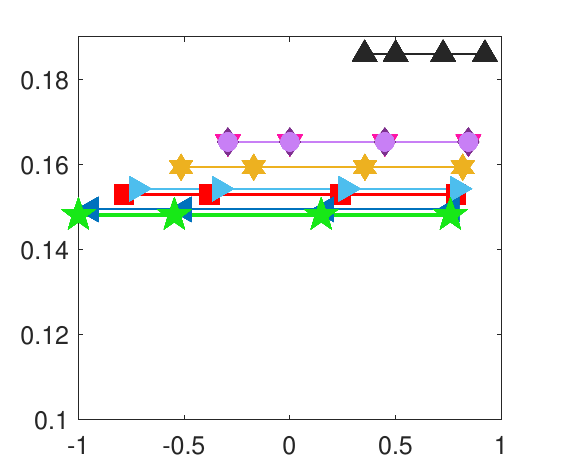} 
      & \includegraphics[width=3.5cm]{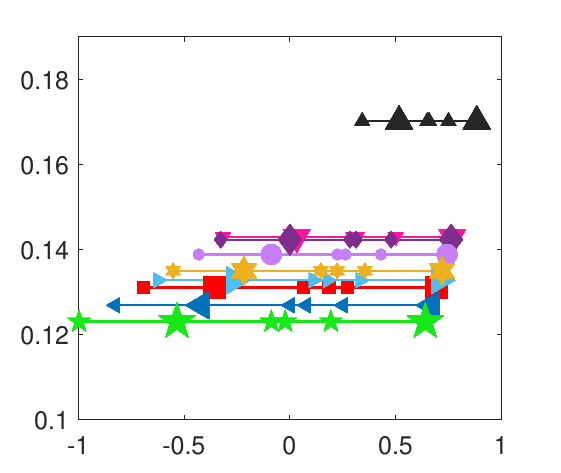} \\
     
      \begin{tabular}{@{}l@{}}
         Acc-DNGD \\ 
         \scriptsize{$E_{\S} = $ worst-case of} \\
         \scriptsize{$f(\xb^K)-f(x^*)$} \\
         \scriptsize{for $K=5$.}
      \end{tabular}
     & \includegraphics[width=3.5cm]{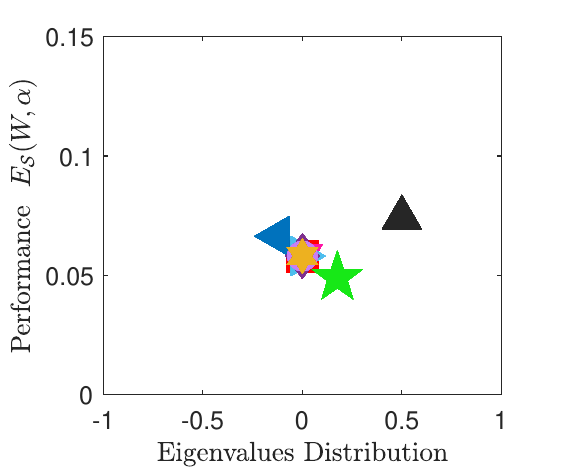} 
     & \includegraphics[width=3.5cm]{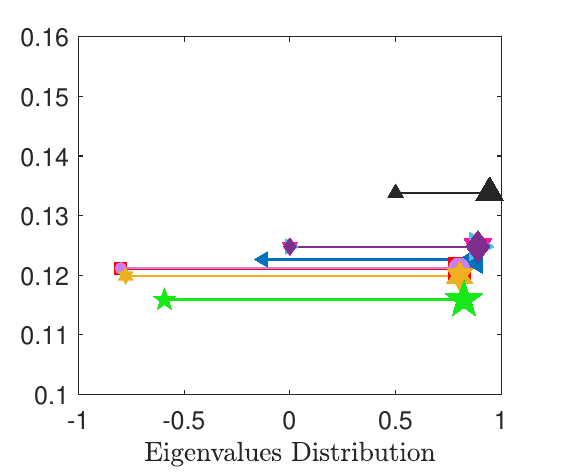} 
     & \includegraphics[width=3.5cm]{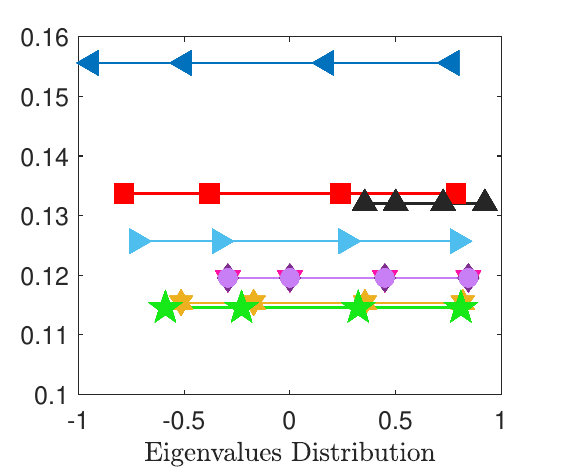} 
     & \includegraphics[width=3.5cm]{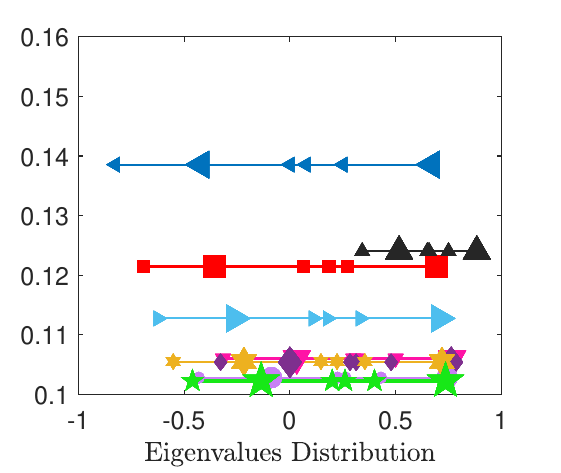} \\
   \end{tabular}
   \caption{These plots show the error criterion $E_{\S}(W,\alpha)$ on the vertical axis, for the optimal averaging matrix $W^*$ in comparison with different averaging matrix heuristics from Section \ref{sec:heuristics}. The plots also show the eigenvalue distribution of the matrices (except $\lam_1=1$) on the horizontal axis. Each marker corresponds to a different eigenvalue, with a size proportional to its multiplicity. To obtain a fair comparison between the averaging matrices, we tune the step-size $\alpha$ of the methods for each of them. Each plot corresponds to a different topology or a different algorithm. The local functions are $\mu$-strongly convex and $L$ smooth. \hfill~ }
   \label{fig:results} 
 \end{figure*}

Using the \emph{pattern search} method from Matlab, we solve problem \eqref{eq:tuning2} for different performance settings $\S$. Fig. \ref{fig:results} compares the optimal weights with the heuristics weights presented in Section \ref{sec:heuristics}. To obtain a fair comparison, we tune the step-size $\alpha$ of the methods, for each heuristic. 
We analyze averaging matrices for the following algorithms, with a constant step-size $\alpha$,
\begin{itemize}
   \item DIGing \cite{DIGing}, see \eqref{eq:DIGing1} and \eqref{eq:DIGing2},
   \item ATC-DIGing \cite{ATC-DIGing},
   \item EXTRA \cite{EXTRA}, with $\tilde{W} = \frac{I+W}{2}$,
   \item Acc-DNGD \cite{AccDNGD}.
\end{itemize}

We have decided to analyze these algorithms because of their significance in the literature.
As explained in \cite{Colla_PEP_dec}, for DIGing and ATC-DIGing, we can compute a bound on the linear convergence rate using \eqref{eq:PEP} by analyzing only $K=1$ iteration of the algorithm, with a performance criterion identical to the initial condition. For the other algorithms (EXTRA and Acc-DNGD), we compute the decrease after $K=5$ steps.
We assume local functions to be $\mu$-strongly convex and $L$-smooth, so that all the above algorithms exhibit a linear convergence.
We consider 4 usual network topologies for $n$ nodes: complete, star, cycle, and grid. Details of these topologies are given in the table header from Fig. \ref{fig:results}.
We mainly focus on the situation with $n=9$ nodes. This small number of nodes gives us a better intuition of the results and also avoids excessive computational load to solve the tuning problem \eqref{eq:tuning2}.
To further reduce this computational load, we have decided to give equal weights to equivalent edges (see Definition \ref{def:eqedges}), so that the actual number of variables in the problem is very low. This corresponds to a restriction on the searching space of \eqref{eq:tuning2}, which may lead to suboptimal solutions if the problem is not convex in $W$.
\begin{definition}[Equivalent edges, \cite{gross2018graph}] \label{def:eqedges} ~\\
   Let $G(\V,\E)$ be a graph and $e_1, e_2 \in \E$ two edges of the graph. If there is an automorphism mapping $e_1$ to $e_2$, then there is an equivalence relation between $e_1$ and $e_2$. The equivalence classes of edges of $G$ are called the edge orbits.  
\end{definition}
For example, the complete graph only has one orbit with all the edges. This is also the case for the star graph and the cycle graph. Therefore, we only have one constant edge weight to optimize for these topologies. The grid graph with 9 nodes, shown in Fig. \ref{fig:results}, has 2 edge orbits: one with the outer edges and the other with the inner edge. For larger grid graphs, each of these orbits divides in two so that we have 4 edge orbits.
We may expect better optimal values for the tuning problem \eqref{eq:tuning2} when allowing all the weights to be different, however, we have observed on small networks ($n=3,4,5$), that this is not the case for the algorithms considered and that we still obtain equal weights for equivalent edges, even when not imposed in advanced.
Moreover, all the weights heuristics presented in Section \ref{sec:heuristics} satisfy this equivalent weight property. \vspace{-2mm} 

\subsection{Weights minimizing $\lam$ (SLEM) are not optimal} \vspace{-2mm}
As shown in Fig. \ref{fig:results}, the averaging matrix $\Wls$ \eqref{eq:minlam2}, minimizing $\lam$ (SLEM) for the given graph $G$, does not provide the best convergence guarantees for the tested decentralized algorithms, even with a step-size $\alpha$ specifically optimized for these weights. This observation is counter-intuitive because small $\lambda$ values correspond to well-connected graphs and because $\Wls$ is optimal for pure consensus steps \eqref{eq:cons}. All known theoretical performance bounds in decentralized optimization depend on $\lam$, and improve when $\lam$ decreases, however, our results show that for a given topology, minimizing $\lam$ does in general not optimize worst-case performance. This suggests that $\lam$ may not be the best characterization of the network performance in distributed optimization methods.
Even for the complete graph, where a reasonable person would guess that $\Wls=\frac{\mathbf{11}^T}{n}$ would be optimal, we can obtain improvement on the performance by choosing $W^*$ with a different unique eigenvalue.

Choosing the optimal averaging matrix $W^*$ (in green in Fig. \ref{fig:results}), instead of $\Wls$ (in red) improves the iterates convergence rate $\rho$ of DIGing and ATC-DIGing by up to 5\% and their corresponding convergence time $\tau = 1/\log(1/p)$ is up to 4 times smaller, when choosing the best possible step-size for each matrix, as we do in Fig. \ref{fig:results}.
For EXTRA and Acc-DNGD, their worst-case performance after 5 iterations is improved by up to 18\%.
We have observed that the optimal weights $W^*$ allows the algorithms to work better with larger step-sizes $\alpha$, which may partly explain their improved performance.

We observe that the sign of the eigenvalues of the averaging matrix and their distribution are related to the resulting performance. Indeed, in many cases (DIGing, ATC-DIGing, Acc-DNGD), the optimal averaging matrices have a smaller range of eigenvalues $|\lam_2-\lam_n|$, with the leftmost (negative) eigenvalue $\lam_n$ significantly larger than that of $\Wls$, while the rightmost (positive) eigenvalue $\lam_2$ is only slightly larger, which can barely be seen on the plots from Fig. \ref{fig:results}. Matrix $\Wls$ always have $\lam_n = - \lam_2$ but by taking advantage of non-symmetric range of eigenvalues, $W^*$ reach smaller $|\lam_2-\lam_n|$.  For EXTRA, things are surprisingly different, the optimal averaging matrix $W^*$ has a larger range of eigenvalues, the smallest of which is close to -1. This should be investigated to see if the phenomenon persists for larger numbers $K$ of iterations.

\subsection{Weights heuristics analysis}
As shown in Fig. \ref{fig:results}, none of the heuristics from Section \ref{sec:heuristics} achieves optimal weights performance. However, some of them come close in certain contexts, while others never seem to perform well. For example, The lazy-Metropolis weights $W_{\mathrm{lazy}-M}$ (in black) do not perform very well in the settings we tested, because it gives averaging matrices with too large positive eigenvalues. Its eigenvalues are concentrated in a smaller range but always have the largest (positive) eigenvalue among all the tested matrices. The standard Metropolis weights $W_{M}$ present better results, while it can still be far from the optimal weights in some cases. Surprisingly, $W_M$ often beats $\Wls$, especially for the grid topology.

For clarity, Fig. \ref{fig:results} omits the performance of the two heuristics for uniform edge weights \eqref{eq:constant-weight} and \eqref{eq:max-degree-weights}. Their respective performance is always worst or equal to those of $\Wls$ and $W_M$.

We observe that the averaging matrices that perform similarly tend to have the same eigenvalue distribution and have actually weights close to each other, even if similar matrices are not guaranteed to be equal. Moreover, when the performance of a heuristic weight approaches the optimal performance, it also approaches the optimal weight values. This would be consistent with the convexity of problem $\eqref{eq:tuning}$ and the uniqueness of its optimal solution, which cannot be guaranteed theoretically at this stage.

No heuristic strictly outperforms the others in all situations. For DIGing, ATC-DIGing, and Acc-DNGD, the best heuristic seems to be $\Wlms$, which performs very well, except for the star graph. By definition \eqref{eq:least_dev}, $\Wlms$ minimizes the steady-state mean square deviation $\delta_{ss}$ of a consensus with additive noise. In distributed optimization, noise is replaced by all kinds of local updates, which take local gradients into account.
In terms of eigenvalues, $\Wlms$ minimize $\delta_{ss} = \sum_{i=2}^n \frac{1}{1-\lam_i(W)^2}$, which drives the matrix to have small eigenvalues in absolute value and penalizes more the large eigenvalues due to the square.
This explains why the eigenvalues of $\Wlms$ often lie in a smaller range around zero, than other matrices.

For EXTRA, the best heuristic seems to be $\Wr$, which selects the weights leading to the smallest (negative) eigenvalues to minimize $\Rtot$, defined in \eqref{eq:Rtot}. This strategy does not appear to be good for other algorithms.  

The observations made throughout this section have also been validated on two random graphs, sampled from the Erdős–Rényi model $\mathcal{G}(n,p)$, one for $n=9$ and $p=0.4$ and the other for $n=30$ and $p=0.2$.

\section{Conclusion}
We showed how to compute the optimal communication weights for a distributed optimization algorithm over an undirected network. Our analysis reveals that the weights minimizing the second-largest eigenvalue modulus (SLEM) of the averaging matrix are suboptimal. While the SLEM characterizes well the convergence rate of a pure consensus protocol, we showed that this is not the best determinant for the performance of a weight matrix in distributed optimization. While other heuristics give better results, the best characterization of the network performance for distributed optimization is still an open question and is probably involving all the eigenvalues of the averaging matrix. 

%and is largely used in convergence analysis of distributed optimization algorithm, does not the best d  

%\blue{Faire une note quelque part sur les matrices et les step-sizes time-varying.}

\section*{Acknowledgements}
We would like to thank Paul-Victor Coulon for his initial work on this topic in his master thesis.

\bibliography{ifacconf}

\end{document}